\documentclass[11pt,a4paper,leqno,oneside]{article}

\usepackage{amsmath}
\usepackage{amsfonts}
\usepackage{amssymb}
\usepackage{amsthm}
\usepackage{makeidx}
\usepackage{graphicx}
\usepackage{pb-diagram}
\usepackage{epstopdf}
\usepackage{fancyhdr}
\usepackage{enumitem}
\usepackage[all]{xy}

\newtheorem{teo}{Theorem}%definizioni che danno luogo ai nomi di teoremi, proposizioni, ecc.
\newtheorem{prop}[teo]{Proposition} %la quadra [teo] indica che la numerazione  subordinata a quella dei teoremi
\newtheorem{lem}[teo]{Lemma}
\newtheorem{cor}[teo]{Corollary}
\theoremstyle{definition}
\newtheorem{defin}[teo]{Definition}

\newtheorem{nt}[teo]{Notation}

\newtheorem{nota}[teo]{Note}

\newtheorem{con}[teo]{Conjecture}
\newtheorem{prob}[teo]{Problem}

\newtheorem{rem}[teo]{Remark}

   {\begin{verse}%
     \small }%
   {\end{verse}}

\begin{document}

\setlength{\parskip}{1ex plus 0.5ex minus 0.2ex}
\begin{center}
$\mathbf{AN\,\,\,ALGEBRAIC\,\,\,CHARACTERIZATION\,\,\,OF\,\,\,HILBERT\,\,\,LATTICES}$
\end{center}
\begin{center}
V.CAPRARO
\end{center}
$\mathbf{Abstract.}$ In this paper we give an algebraic
characterization of the projections lattice of $M_n(\mathbb C)$ and
we extend it to the case of $B(H)$, with $H$ separable Hilbert
space.

\tableofcontents

\section{Introduction}
In this paper we give an algebraic characterization of the
projections lattice of $M_n(\mathbb C)$ and we extend it to the case
of $B(H)$, with $H$ separable Hilbert space. Such a characterization
was founded by several authors following, above all, a topological
way (see \cite{So},\cite{Pr},\cite{Dv},\cite{Zi}). The reason of
this interest come from the works G.W. Mackey (\cite{Ma}) and
Birkhoff-von Neumann (\cite{B-vN}), in which they axiomatize the
quantum logic with the lattice of closed subspaces (or equivalently
the lattice of projections) of a separable Hilbert space. But we
have in mind to use a characterization of Hilbert lattices in order
to attach the Connes' embedding conjecture (\cite{Co}) via lattice
theory and using a very famous theorem of Kirchberg (\cite{Ki}). For
this reason we need a different characterization of Hilbert
lattices. It should be based only on algebraic properties and it
should describe in the most precise way what happens in the finite
dimensional case, i.e. in the case of the $n\times n$ matrices on
the complex field. The nearest approach we meet is strangely the
first one: during the years 1935-38 J. von Neumann found an
axiomatization for the projections lattice of a finite factor, but
he used some axioms that we will not use. In particular we don't
assume the existence of a transition probability, cutting off at
least ten of the eighteen von Neumann's axioms. Obviously our
construction loses of generality: it will be valid only in the case
of finite factor of type I, i.e. $M_n(\mathbb C)$. Successively we
are able to extend our construction to the separable case. At last
we describe a possible second theorem of correspondence that should
refine the first one arriving to the ''minimal axiomatization'' of
the projections lattice of $M_n(\mathbb C)$, in which we are able to
cut off other axioms.
\\\\
So, our approach is similar to the von Neumann's one. This is why we
want to compare these two approaches in this preliminary section.
\\\\
Let us recall the following
\begin{defin}{\bf (von Neumann, \cite{vN-H2})}\label{cg}
A continuous geometry with transition probability is a system
$(L,\leq,^{\perp},P)$ verifying the following axioms
\begin{enumerate}
\item $\leq$ is a partial ordering on $L$.
\item Each subset of $L$ admits greatest lower bound and least
upper bound. We set $l\vee l'=sup\{l,l'\}$, $l\wedge l'=inf\{l,l'\},
0=inf(L), 1=sup(L)$.
\item $\vee$ and $\wedge$ are continuous in the following sense
\begin{enumerate}
\item If $\{l_i\}$ is an increasing net in $L$ and $l\in L$, then
$$
\bigvee_i(l\wedge l_i)=l\wedge\bigvee_il_i
$$
\item If $\{l_i\}$ is a decreasing net in $L$ and $l\in L$,
then
$$
\bigwedge_i(l\vee l_i)=l\vee\bigwedge_il_i
$$
\end{enumerate}
\item The modular property holds
$$
l\leq l''\Rightarrow(l\vee l')\wedge l''=l\vee(l'\wedge
l'')\,\,\,\forall l'\in L
$$
\item $l\rightarrow l^{\perp}$ is an involutory anti-automorphism of
$L$, i.e. $l^{\perp\perp}=l$ and $l\leq l'$ implies $l'^{\perp}\leq
l^{\perp}$.
\item $l\leq l^{\perp}$ if and only if $l=0$.
\item We say that $l,l'$ are inverse if $l\vee l'=1$ and
$l\wedge l'=0$. This axiom requires that if $l,l'$ and $l,l''$ are
inverse, it cannot be $l'<l''$.
\item $P:L\setminus\{0\}\times L\rightarrow[0,1]$
\item $P(l,l')=1$ if and only if $l\leq l'$.
\item If $l'\leq l''^{\perp}$ and $l\neq0$, then $P(l,l'\vee
l'')=P(l,l')+P(l,l'')$.
\item We say that a sequence $\{l_n\}$ is convergent to $l\in L$ if
$P(l',l_n)\rightarrow P(l',l)$, for all $l'\in L\setminus\{0\}$.
This axiom requires that each increasing sequence is convergent.
\item Let $\{l_n\}\subseteq L$. If $P(1,l_n)=constant\neq0$,
$P(l_n,l_m)\rightarrow1$ for $n,m\rightarrow\infty$ and $P(l,l_n)$
is convergent for each $l\neq0$, then $\{l_n\}$ is convergent.
\item For each $\varepsilon>0$ there exists a
$\delta=\delta_{\varepsilon}>0$ such that $P(1,l')\leq\delta$
implies $P(1,l\vee l')\leq P(1,l)+\varepsilon$ for all $l\in L$.
\item If $l\in L$ satisfies $l'=(l'\wedge l)\vee(l'\wedge
l^{\perp})$ for all $l'\in L$, then $l=0$ or $l=1$. This property is
called ''irreducibility''.
\item If $T$ is an $(L,\leq,\perp)$-automorphism, then it is
also a $P$-automorphism, i.e. $P(Tl,Tl')=P(l,l')$, for all $l,l'\in
L, l\neq0$.
\item If $P(1,l)<P(1,l')$ then there exists an $(L,\leq,\perp)$-automorphism $T$ such that $Tl<l'$.
\item If $P(1,l)=P(1,l')$ then there exists an $(L,\leq,\perp)$-automorphism $T$ such that $Tl=l'$ and $Tl''=l''$ for each
$l''\leq(l\vee l')^{\perp}$.
\item There exist $l_1<l_2<l_3<l_4$.
\end{enumerate}
\end{defin}
In \cite{vN-H2} von Neumann proved that a continuous geometry with
transition probability is isomorphic to the projection lattice of a
finite factor (type $I_n$ or $II_1$ depends by the values of the
dimension
function, whose definition we will recall later).\\
In this paper we obtain the same result in the case $I_n$, but
without using the transition probability. So we are able to cut off
the axioms: 8, 9, 10, 11, 12, 13, 15, 16, 17. Moreover we will not
use axiom 18 and, apparently, axioms 7 (apparently means that one
should study the construction of dimension function deeper in order
to see what axioms von Neumann used). On the other hand, we  have to
add two obvious axioms: $L$ must have the same cardinality of
$\mathbb R$ and it must contain at least a minimal element.
Actually, we will give this last axiom in a (we hope!) more elegant
way, generalizing the well-known property of von Neumann algebras:
the center of the reduced algebra is the reduced algebra of the
center. But we will use this property only to prove that there exist
at least one minimal element.
\section{Orthocomplemented lattices}
In this section we recall the original von Neumann's definition of
complete lattice (cfr. \cite{vN-H1}).
\begin{defin}
A complete lattice is a system $(L,\leq)$ verifying the following
axioms:
\begin{enumerate}
\item $\leq$ is a partial ordering on $L$.
\item Each subset of $L$ admits a least upper bound and a greatest lower
bound. We set $sup(L)=1$, $inf(L)=0$, $inf(l,l')=l\wedge l'$ and
$sup(l,l')=l\vee l'$.
\item $\vee$ and $\wedge$ are continuous in the following sense
\begin{enumerate}
\item If $\{l_i\}$ is an increasing net in $L$ and $l\in L$, then
$$
\bigvee_i(l\wedge l_i)=l\wedge\bigvee_il_i
$$
\item If $\{l_i\}$ is a decreasing net in $L$ and $l\in L$,
then
$$
\bigwedge_i(l\vee l_i)=l\vee\bigwedge_il_i
$$
\end{enumerate}
\item The modular property holds
$$
l\leq l''\Rightarrow(l\vee l')\wedge l''=l\vee(l'\wedge
l'')\,\,\,\forall l'\in L
$$
\end{enumerate}
\end{defin}
\begin{rem}
More recently one calls lattice only a partially ordered set in
which each pair of elements admits sup and inf. In particular a
lattice is called modular if it verifies the modular property. We
work only on modular lattices (unless a little example in the next
section) thus we preferred to follow the original von Neumann's
definition. On the other hand, it is just the modular property that
make false our characterization in infinite dimension: it is proved
in \cite{B-vN} that the projections lattice of $B(H)$, with $H$
infinite-dimensional, is not modular. (see also \cite{Fu}). It might
be interesting to see in details what happens without this axiom. We
think that our construction can be generalized (maybe generalizing
the dimension function), since it is based above all on minimal
elements/projections. We will discuss this idea in the last section.
\end{rem}
\begin{defin}
An orthocomplemented lattice is a pair $(L,\perp)$, where $L$ is a
complete lattice and $^\perp:L\rightarrow L$ is an involution that
satisfies
\begin{enumerate}
\item $l\vee l^{\perp}=1$
\item $l\wedge l^{\perp}=0$
\item $l\leq l'$ implies $l'^{\perp}\leq l^{\perp}$
\end{enumerate}
We write also $1-l$ instead of $l^{\perp}$.
\end{defin}
\begin{nota}
In the definition of continuous lattice with transition probability,
von Neumann does not require properties 1. and 2. But one can easily
prove that they follow by the von Neumann's sixth axiom and by the
following well-known properties.
\end{nota}
\begin{enumerate}
\item $(l\vee l')^{\perp}=l^{\perp}\wedge l'^{\perp}$
\item $(l\wedge l')^{\perp}=l^{\perp}\vee l'^{\perp}$
\end{enumerate}
Let $L$ be an orthocomplemented lattice.
\begin{defin}
Let $l\in L$. We set $\perp(l)=\{l'\in L:l'\leq1-l\}$. Elements
belonging into $\perp(l)$ are said to be orthogonal to $l$.
\end{defin}
\begin{defin}
Let $l,l'\in L$. We said that $l$ commutes with $l'$ if
$$
l=(l\wedge l')\vee(l\wedge l'^{\perp})
$$
$c(l)$ will denotes the set of elements which commute with $l$.
\end{defin}
It is well-known that
\begin{enumerate}
\item $\{0,1,l\}\subseteq c(l)$
\item $l'\in c(l)\Leftrightarrow l\in c(l')$
\end{enumerate}
Thus the relation of commutation is symmetric and we can say that
''l,l' commute''.
\begin{defin}
We set
$$
C(L)=\bigcap_{l\in L}c(l)
$$
Elements belonging into $C(L)$ are called central. Certainly
$C(L)\supseteq\{0,1\}$. L is called factorial if $C(L)=\{0,1\}$,
abelian if $C(L)=L$.
\end{defin}
\begin{nota}
Factoriality corresponds to irreducibility (14th von Neumann's
axiom).
\end{nota}
\begin{lem}\label{compl}
Let $l\in L$. Then $\perp(l)\subseteq c(l)$
\begin{proof}
Let $l'\leq1-l$ and thus $l\leq1-l'$. We have
$$
[l\wedge(1-l')]\vee(l\wedge l')=l\vee0=l
$$
thus $l'\in c(l)$.
\end{proof}
\end{lem}
\begin{prop}\label{distributiva}
Let $l_3\in c(l_1)\cap c(l_2)$, then
$$
(l_1\vee l_2)\wedge l_3=(l_1\wedge l_3)\vee(l_2\wedge l_3)
$$
\begin{proof}
Applying the modular law (sometimes!), one has
$$
(l_1\vee l_2)\wedge l_3=[l_1\vee(l_2\wedge l_3)\vee(l_2\wedge
l_3^{\perp})]\wedge l_3=
$$
$$
=(l_2\wedge l_3)\vee[(l_1\vee(l_2\wedge l_3^{\perp}))\wedge l_3]=
$$
$$
=(l_2\wedge l_3)\vee[((l_1\wedge l_3)\vee(l_1\wedge
l_3^{\perp})\vee(l_2\wedge l_3^{\perp}))\wedge l_3]=
$$
$$
=(l_2\wedge l_3)\vee[(l_1\wedge l_3)\vee(((l_1\wedge
l_3^{\perp})\vee(l_2\wedge l_3^{\perp}))\wedge l_3)]=
$$
$$
=(l_2\wedge l_3)\vee(l_1\wedge l_3)
$$
\end{proof}
\end{prop}
\begin{cor}\label{abeliano}
An othocomplemented lattice $L$ is abelian if and only if for each
$l,l',l''\in L$ the following holds
$$
(l\vee l')\wedge l''=(l\wedge l'')\vee(l'\wedge l'')
$$
Thus abelian lattices coincide with those are classically called
distributive lattices or Boole algebras.
\begin{proof}
If $L$ is abelian, one can apply Prop. \ref{distributiva}.
Conversely, one can apply the formula with $l'=1-l$.
\end{proof}
\end{cor}
\begin{prop}
Let $l,l'\in\perp(l)$. If $l\vee l'=l\vee l''$, then $l'=l''$.
\begin{proof}
By Lemma \ref{compl} we have $l\in c(l')\cap c(l'')$. Thus, using
Prop. \ref{distributiva},
$$
l''=(l\wedge l'')\vee(l''\wedge l'')=(l\vee l'')\wedge l''=
$$
$$
=(l\vee l')\wedge l''=(l\wedge l'')\vee(l'\wedge l'')=l'\wedge l''
$$
Whence $l''=l'\wedge l''$ and consequently $l''\leq l'$. Changing
$l'$ and $l''$ we can find the reverse inequality.
\end{proof}
\end{prop}
\begin{rem}
Let $l\in L$, we consider the complete lattice $L'=L\wedge
l=\{l'\wedge l, l'\in L\}$. This can be orthocomplemented in a
natural way setting (for each $l'\leq l$)
$$
l-l'=(1-l')\wedge l
$$
Notice that one have to use the modular property to prove that
$l'\vee(l-l')=l$. Indeed
$$
l'\vee(l-l')=l'\vee[(1-l')\wedge l]=[l'\vee(1-l')]\wedge l=1\wedge
l=l
$$
\end{rem}
\begin{defin}
An element $l\in L$ is called abelian if $L\wedge l$ is an abelian
lattice. Let $Ab(L)$ be the set of abelian element in $L$.
\end{defin}
\begin{rem}
Obviously a Boole algebra can be characterized as a lattice for
which $Ab(L)=L$.
\end{rem}
\begin{defin}
A factorial lattice is said ''of type I'' if it admits at least a
non-zero abelian element.
\end{defin}
\begin{rem}
By definition of orthocomplement in $L\wedge l$, it follows that
$C(L)\wedge l\subseteq C(L\wedge l)$. The reverse inclusion holds
for Boole algebras and it could hold for other important classes of
orthocomplemented lattices.
\begin{defin}
An $R$-lattice is an orthocomplemented lattice in which the
following property of restriction of the central elements holds
$$
C(L\wedge l)=C(L)\wedge l\,\,\,\,\,\forall l\in L
$$
\end{defin}
We repeat that this property generalized the well-known property of
von Neumann algebras: the center of the reduced algebra is the
reduced algebra of the center.
\begin{prob}
Do there exist orthocomplemented lattices which are not R-lattices?
also among the factorial ones?
\end{prob}
\end{rem}
We conclude this preliminary section giving a nice characterization
of the orthogonal complement. We recall the following
\begin{defin}\label{inv}
Let $l\in L$. An inverse of $l$ is an element $l'\in L$ such that
$l'\vee l=1$ e $l\wedge l=0$. $inv(l)$ will denote the set of the
inverse of $l$.
\end{defin}
\begin{prop}\label{ort}
Let $l\in L$. One has
$$
c(l)\cap inv(l)=\{1-l\}
$$
\begin{proof}
Certainly $1-l\in c(l)\cap inv(l)$. Now let $x\in c(l)\cap inv(l)$,
then
$$
x=(x\wedge l)\vee(x\wedge(1-l))=x\wedge(1-l)
$$
thus $x\leq1-l$. Now we assume that $x\neq1-l$, then $1-l-x$ is
defined and not zero. Now we observe that if $l,l'$ are two elements
in $L$ and $l'\leq l$, we defined $l-l'=l\wedge(1-l')$. Now we have
$1-l'\leq l$ and thus
$$
l-(1-l')=l\wedge(1-(1-l'))=l\wedge l'
$$
Applying this observation with $l=1-l$ and $l'=1-x$ we have
$$
1-l\vee x=(1-l)\wedge(1-x)=1-l-x\neq0
$$
ant thus the absurd: $l\vee x\neq1$.
\end{proof}
\end{prop}
\begin{defin}\label{ind}
$l,l'\in L$ are said to be independent if $l\wedge l'=0$.
\end{defin}
\begin{rem}
To prove the implication $x\in c(l)\cap inv(l)\Rightarrow x\leq1-l$
we used only the independence between $l$ and $l'$. Thus we can
interpret $1-l$ as the sup of elements which are independents and
commuting with $l$ .
\end{rem}
\section{Some examples}\label{esem}
In this section we give some examples of othocomplemented lattices.
\begin{enumerate}
\item For each set $X$, the set of its subsets $P(X)$ is a Boolean algebra. It is
well-known (see \cite{Bi}) that each finite Boolean algebra can be
obtained in this way. The Stone representation theorem (see
\cite{St} or \cite{Jo}) give the following deep characterization:
each Boolean algebra is isomorphic to the set of clopen subsets
(subsets which are open and closed) of an Hausdorff extremely
disconnected topological space.
\item The projections lattice of a von Neumann algebra is always an
infinite orthocomplemented lattice. It is abelian (resp. factorial)
if and only if the von Neumann algebra is abelian (resp. a factor).
\item Let $n\in\mathbb N$. If $n$ is odd, there are no orthocomplemented lattices
with $n$ elements. On the other hand, if $n$ is even, and equal to
$2m$, it is possible to construct a factorial lattice with $n$
elements. Indeed we set
$$
L_m=\{0,l_1,1-l_1,l_2,1-l_2,...l_m,1-l_m\}
$$
with the conditions $l\wedge l'=0, l\vee l'=1$ for each $l\neq l'$.
This lattice can be represented in the following way
$$
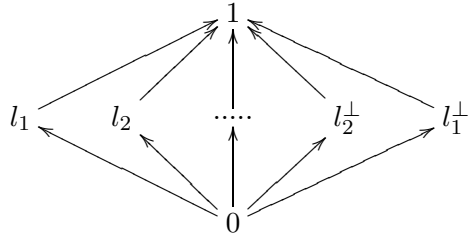
\xymatrix{& & 1 & \\
l_1\ar[urr] & l_2\ar[ur] & ..... \ar[u] & l_2^{\perp}\ar[ul] &
l_1^{\perp}\ar[ull]\\
& & 0\ar[ull] \ar[ul] \ar[u] \ar[ur] \ar[urr] }
$$
Chevalier showed that with these lattices, together with the
distributive ones, one can construct each finite orthocomplemented
lattice. Indeed one can define (in a obvious way) the algebraic
product of a finite family of orthocomplemented lattices with the
lexicographic order and setting
$(l_1,...l_n)^{\perp}=(l_1^{\perp},...l_n^{\perp})$. Chevalier
proved that each finite orthocomplemented lattice is product of
$\{0,1\}^n$ with a finite family of type $L_m$ lattices. (see also
\cite{Sv}).
\item More recently lattices which are not modular are studied.
The simplest of those is the pentagon $N_5$
$$
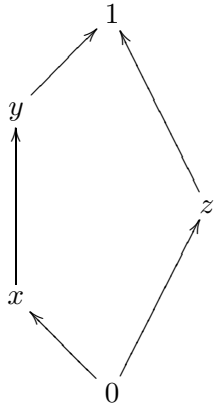
\xymatrix { & 1 \\
y\ar[ur] & &\\
& & z\ar[uul]\\
x\ar[uu] & &\\
&0\ar[ul] \ar[uur] }
$$
It does not satisfies the modular property. Indeed
$$
(x\vee z)\wedge y=y\neq x=x\vee(z\wedge y)
$$
In \cite{G} it is proved that a lattice verifies the modular
property if and only if it does not contain sublattices isomorphic
to $N_5$.
$$
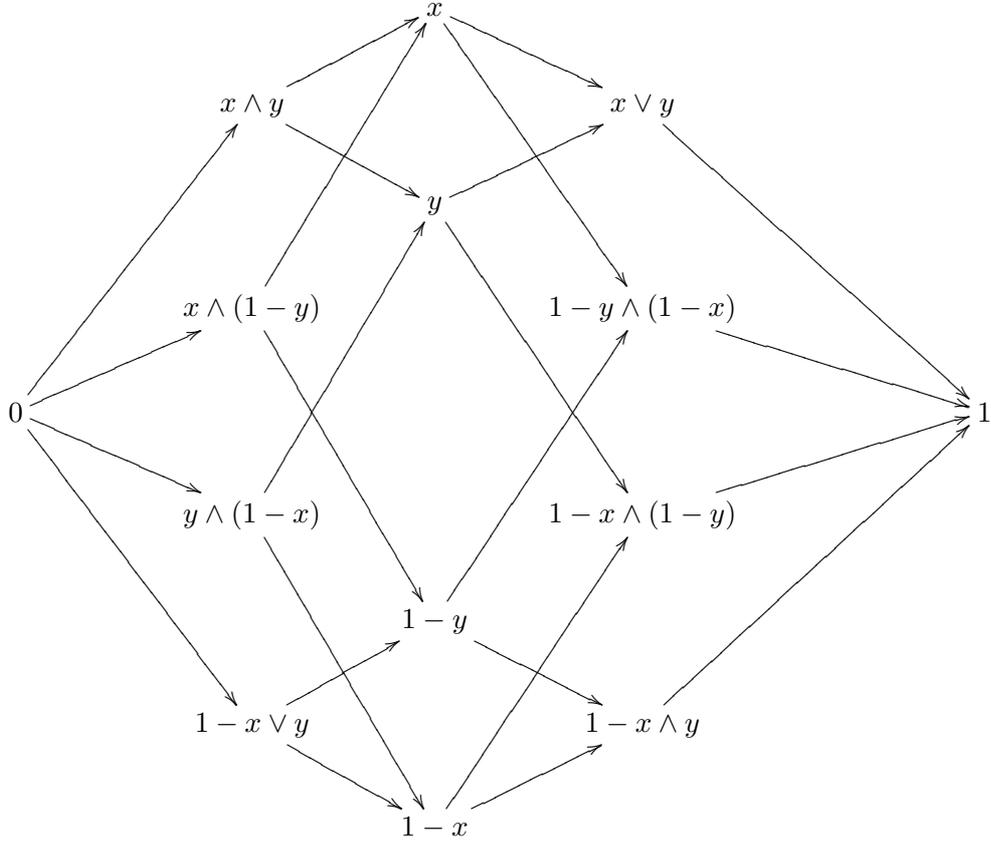
\xymatrix{ & & & x \ar[dr] \ar[dddr] &\\
& &x\wedge y \ar[ur] \ar[dr] & & x\vee y \ar[dddrrr]\\
& & & y \ar[ur] \ar[dddr]& \\
& & x\wedge(1-y)\ar[uuur]\ar[dddr] & & 1-y\wedge(1-x)\ar[drrr]\\
0 \ar[uuurr] \ar[urr] \ar[drr] \ar[dddrr]& & & & & & & 1\\
& & y\wedge(1-x)\ar[uuur] \ar[dddr] & & 1-x\wedge(1-y)\ar[urrr]\\
& & & 1-y \ar[uuur]\ar[dr]& \\
& &1-x\vee y \ar[ur] \ar[dr] & & 1-x\wedge y\ar[uuurrr]\\
 & & & 1-x \ar[ur] \ar[uuur] &}
$$
Notice that $N_5$ is not an orthocomplemented lattice. To construct
an orthocomplemented lattice which does not verify the modular
property one can consider the previous one. Indeed it is clearly an
orthocomplemented lattice, in the sense that it satisfies the first
three property of the definition of lattice and all those of the
definition of the orthocomplementation. But it is not modular (you
find with $\{x,1-x,x\vee y\}$). Moreover, a direct calculation shows
that $C(L)=L$. On the other hand it is not distributive (otherwise
it should be modular!). Notice that, in order to obtain the
equivalence in Cor.\ref{abeliano}, we used the modular property to
prove Prop. \ref{distributiva}. More generally the following
proposition holds
\begin{prop}
Let $L$ be an orthocomplemented non-modular lattice. One considers
the following three properties
\begin{enumerate}
\item $L$ is modular.
\item $L$ is abelian.
\item $L$ is distributive.
\end{enumerate}
Then $c)\Leftrightarrow a)\wedge b)$. That is: the third property
holds if and only if both the other two hold.
\end{prop}
\begin{prob}
What is the smallest example of orthocomplemented lattice which is
neither abelian nor factorial?
\end{prob}
Notice that the previous one is not the smallest example of
orthocomplemented lattice which is abelian without being
distributive. The hexagon is a much more trivial example!
\end{enumerate}
\section{Equivalence relations on orthocomplemented lattices}
In \cite{vN-H1} von Neumann showed that a powerful mean to know the
structure of an orthocomplemented lattice is given by its
equivalence relations. In the next three sections we want to run
again the street opened by von Neumann, abstracting some concepts.
\begin{defin}
Let $L$ be a lattice and $\sim$ an equivalence relation on $L$. We
say that $l\in L$ is $\sim$dominated by $l'\in L$ (and we write
$l\leq_{\sim}l'$) if and only if there exists $l''\leq l'$ such that
$l''\sim l'$. We write $l<_{\sim}l'$ if $l''<l'$.
\end{defin}
Obviously one is interested to equivalence relation that are
compatible with the ordering in some sense.
\begin{defin}
An equivalence relation on an orthocomplemented lattice is called
regular if
\begin{enumerate}
\item $l\sim0\Leftrightarrow l=0$
\item $l\geq l'$ and $l\leq_{\sim}l'$ imply $l\sim l'$
(order compatibility)
\item For each $l,l'\in L$ one and only one of the followings hold
$l\sim l'$, $l\leq_{\sim}l'$, $l'\leq_{\sim}l$.
\item If $\{l_i\}$ and $\{l'_i\}$ are two families of mutually
orthogonal elements and they are such that $l_i\sim l_i'$ for each
$i$, then
$$
\bigvee_il_i\sim\bigvee_il_i'
$$
\item Conditions $l'\leq l, l'\sim l$ imply $l'=l$ (finiteness
property).
\end{enumerate}
\end{defin}
\begin{defin}\label{fattore}
A complete lattice is called irreducible if only $0$ and $1$ have
unique inverse (see Def.\ref{inv})
\end{defin}
\begin{rem}\label{irr}
We recall that irreducibility corresponds (in orthocomplemented
lattices) to factoriality (see. \cite{vN-H2}, pg. 14).
\end{rem}
We can join some von Neumann's results to obtain the following
\begin{teo}{\bf (von Neumann)}\label{vn}
Let $L$ be a complete irreducible lattice. Setting $l\sim l'$ if and
only if they are a common inverse (see Def.\ref{inv}), then $\sim$
is a regular equivalence relation on $L$.
\begin{proof}
Properties 2,3,4 e 5 are already proved by von Neumann in
\cite{vN-H1}. The first one is clear: if $l\sim0$, then $l,0$ have a
common inverse. But the unique inverse of $0$ is $1$ and thus
$l=1\wedge l=0$.
\end{proof}
\end{teo}
\begin{defin}
We say that two elements are perspective if they are a common
inverse.
\end{defin}
Thus each irreducible complete lattice (and in particular each
orthocomplemented factorial lattice) admits a regular equivalence
relation. In this case one can ask if this relation is also
compatible with the orthogonality in some sense. Property 4. is
already a compatibility with the orthogonality, but we will prove
something more: there exists substantially only one regular relation
(which coincide with the perspective relation, i.e.
with the Murray-von Neumann relation on projections) and thus it verifies the parallelogram law: $l\vee l'-l\sim l'-l\wedge l'$.\\
\begin{prop}\label{ordine}
Let $\sim$ be a regular equivalence relation on $L$. If $l\sim l'$
and $l_1\leq l$, then $l_1\leq_{\sim}l'$.
\begin{proof}
By definition one has $l_1\leq_{\sim}l'$ or $l'\leq_{\sim}l_1$. We
assume the last one and we prove that actually the first one holds
in the form: $l_1\sim l'$. Indeed $l\sim l'\leq_{\sim}l_1$ and thus,
by transitivity, $l\leq_{\sim}l_1$. On the other hand, by
hypothesis, $l_1\leq l$ and thus, by order compatibility, it follows
that $l_1\sim l$ and thus $l_1\sim l'$.
\end{proof}
\end{prop}
\section{Minimal elements in factorial lattices}
Let $L$ be a factorial $R$-lattice and $\sim$ be a regular
equivalence relation on $L$. We recall the property R: $C(L)\wedge
l=C(L\wedge l)$ for each $l\in L$. We don't know if there exists
factorial lattices which do not verify it.
\begin{defin}
A non-zero element $l\in L$ is called minimal if for each $l'\in L$
one and only one of the following holds: $l\wedge l'=0$, or $l\wedge
l'=l$. Let $Min(L)$ denote the set of minimal elements of $L$ and
with $Min(l)$ the set of minimal elements of $L$ which are less then
$l$.
\end{defin}
\begin{nota}
In literature minimal elements are also called atoms. We prefer to
call them minimal for keeping a sense of continuously with the
theory of Operator Algebras.
\end{nota}
\begin{cor}\label{due}
If $l$ is minimal and $l'\sim l$, than also $l'$ is minimal.
\begin{proof}
Let $l_1'\leq l'$ and thus (by Prop.\ref{ordine}) $l_1'\leq_{\sim}l$
and consequently $l_1'\sim0$, or $l_1'\sim l$. The properties of
$\sim$ then guarantee $l_1'=0$, or, by finiteness, $l_1'=l'$.
\end{proof}
\end{cor}
Here is one of the main result to develop our approach.
\begin{prop}\label{mini}
Let $L$ be a type I factorial lattice. Each element of $L$ contains
a minimal element.
\begin{proof}
Let $l\neq0$, $l\in Ab(L)$. Factoriality and property R imply
$\{0,l\}=C(L)\wedge l=C(L\wedge l)$; and thus $L\wedge l$ is still
factorial. Since it is also abelian, it must be $L\wedge l=\{0,l\}$.
Thus $l$ is minimal. Now let $l'\in L$. It must be $l'\leq_{\sim}l$,
or $l\leq_{\sim}l'$. In the first case $l'\sim l_1\leq l$. Thus
$l'\sim0$ (and therefore $l'=0$), or $l'\sim l$ (and therefore $l'$
is minimal, by Cor.\ref{due}). In the second one $l'$ contains an
element equivalent to $l$. This element must be minimal, being $l$
minimal (still by Cor.\ref{due}).
\end{proof}
\end{prop}
\begin{cor}\label{rife}
Let $L$ be a type I factorial lattice and $\{l_i\}$ a family of
elements which are maximal in the properties to be minimal and
mutually orthogonal. Then $\bigvee_il_i=1$.
\begin{proof}
If $l=1-\bigvee_il_i\neq0$, we can find (see Prop.\ref{mini}) a
minimal element $l_1\leq l$, contradicting maximality.
\end{proof}
\end{cor}
\begin{cor}\label{minimale}
Let $L$ be a type I factorial lattice. Each element of $L$ is sup of
a family of minimal and mutually orthogonal elements.
\end{cor}
\section{The dimension function}
In this section we give a sketch of the construction of the
dimension function in an irreducible complete lattice.
\begin{defin}
A system $\{l_i\}_{i\in I}$ of elements in $L$ is called independent
if for each partition $\{J,K\}$ of $I$ one has
$$
\bigvee_{i\in J}l_i\wedge\bigvee_{i\in K}l_i=0
$$
\end{defin}
\begin{nota}
In this section one can change the word ''independent'' with
''orthogonal''. Th.\ref{vn} is still valid.
\end{nota}
\begin{defin}
Let $L$ be a complete lattice and $\sim$ a regular equivalence
relation on $L$. A dimension function $\sim$-compatible is a map
$D:L\rightarrow[0,1]$ such that
\begin{enumerate}
\item $D(0)=0, D(1)=1$
\item $D(l\vee l')+D(l\wedge l')=D(l)+D(l')$
\item $D(l)=D(l')\Leftrightarrow l\sim l'$
\item $D(l)\leq D(l')\Leftrightarrow l\leq_{\sim}l'$
\item If $\{l_i\}$ is a finite or countable independent system, then
$$
D(\bigvee l_i)=\sum D(l_i)
$$
\end{enumerate}
\end{defin}
Von Neumann proved in \cite{vN-H1} that choosing $\sim$ as the
perspective equivalence relation, then there exists a unique
dimension function $\sim$-compatible. Moreover the image $D(L)$ can
be only one among the following sets
\begin{enumerate}
\item $\Delta_n=\{0,\frac{1}{n},\frac{2}{n}, ... \frac{n-1}{n},1\}$
\item $\Delta_{\infty}=[0,1]$
\end{enumerate}
Here we give a sketch of the construction of the dimension function.
The complete construction can be founded in the first five chapters
of \cite{vN-H1}. Let us precise that this construction also works in
our case (factorial R-lattices), since factoriality and
irreducibility
coincide for orthocomplemented lattices.\\\\
(I step)\\
Since we want define a class function (i.e. constant on each
equivalence class) which describes the order of the equivalence
class through the order of their values, the first step consists of
defining some operations between equivalence classes. More
precisely, we set $\mathcal L=L/\sim$ and we denote $A,B,C...$ the
elements of $\mathcal L$. Von Neumann himself proved that if there
exist $a\in A,b\in B$ such that $a\wedge b=0$ then one can
well-define a unique class $A\vee B$ (that is the class containing
$a\vee b$). Analogously , if there exist $a\in A, b\in B$ such that
$a\geq b$, one can well-define a unique class $A-B$ (that is the class containing $a-b$).\\\\
(II step)\\
After proving some algebraic properties of the operations between
classes, von Neumann gave the following crucial
\begin{defin}
Let $A_0$ be the class containing $0\in L$ and $A\in\mathcal L$. We
set $0A=A_0$. Now, supposing defined $(n-1)A$ and that there exists
$(n-1)A\vee A$, than we define $n A=(n-1)A\vee A$. Otherwise $nA$ is
not defined.
\end{defin}
In order to understand the idea of this definition, it is better to
think at factorial R-lattices of type I. In this case we know that
each element is sup of a family of minimal elements. Now we assume
the following fact (that we will prove later)
\begin{teo}\label{dimensione}
Let $L$ be a factorial R-lattice of type $I$ and $\{l_i\}_{i\in
I},\{l'_j\}_{j\in J}$ two families of mutually orthogonal and
minimal elements. Then
\begin{enumerate}
\item $|I|$ and $|J|$ are finite numbers.
\item $\bigvee_il_i\sim\bigvee_jl'_j$ if and only if $|I|=|J|$.
\end{enumerate}
\end{teo}
Let $A\in\mathcal L$. Thanks to this theorem we can define $m_1$ as
the number of minimal elements appearing in some decomposition of
$1\in L$ and $m_A$ as the number of minimal elements appearing in
some decomposition of $a\in A$. Let $A_{min}$ be the class
containing each minimal element of $L$ and let $n_A$ be the greatest
integer for which $nA$ is defined. In this case von Neumann's idea
is to determine $m_1$ as $n_{A_{min}}$ (indeed $nA_{min}\vee
A_{min}$ is defined until one takes a maximal family of minimal and
mutually orthogonal elements), then to determine
$n_A=n_{A_{min}}-m_A$ and lastly to define the dimension of the
class $A$ as $m_A/n_{A_{min}}$. One can follow this idea also in
more general cases: the key is to understand what means, given two
classes $A$ and $B$: $n$ is the greatest integer for which $nA$ is
into $B$.
\begin{defin}
Let $A,B\in\mathcal L$. We set $A<B$ if and only if there exist
$a\in A,b\in B$ such that $a<_{\sim}b$.
\end{defin}
Now we can enunciate one of the most important theorem of this
theory.
\begin{teo}
Let $A_0\neq A,B\in\mathcal L$. There exist a unique pair $(n,B_1)$,
where $n\geq0$ is an integer and $B_1<A$ is a class such that
$$
B=nA+B_1
$$
We set $[B:A]=n$.
\end{teo}
Now the construction would be concluded in the case of factorial
lattices of type I. In the general case in which we have not minimal
elements, we need a further step.\\\\
(III step)
\begin{defin}
A class $A_0\neq A\in\mathcal L$ is said to be minimal if there no
exists $B\neq A_0$ such that $B<A$.
\end{defin}
The case in which $A$ is not minimal is optimally described by von
Neumann with the following
\begin{teo}\label{minimal}
If $A$ is not minimal, then there exist $B\neq A_0$ such that
$2B\leq A$.
\end{teo}
\begin{defin}
A minimal sequence $\{A_n\}$ of elements $\neq A_0$ is one
containing but one element $B$ which is minimal, or containing a
denumerable infinitude of elements such that $2A_{n+1}\leq A_n$.
\end{defin}
Using th.\ref{minimal} we have the existence of minimal sequence.
The following theorem allows to define the dimension function.
\begin{teo}
Let $A_0\neq A\in\mathcal L$. Then the following limit exists and it
is finite and positive
$$
(A:A_1)=lim_{i\rightarrow\infty}\frac{[A:A_i]}{[A_1:A_i]}
$$
where $A_1$ denotes the class in which $1\in L$ belong. (If
$\{A_i\}$ consists of one minimal element $B$, we mean by
$lim_{i\rightarrow\infty}$ the value at $A_i=B$.).
\end{teo}
(IV step)
\begin{defin}
Let $l\in L$ and $A_l$ be the class containing $l$. We set
$D(l)=(A_l:A_1)$ if $l\neq0$. Otherwise $D(0)=0$.
\end{defin}
The dimension function verifies the following properties.
\begin{enumerate}
\item $D(l)\in[0,1]$ for each $l\in L$. $D(0)=0, D(1)=1$
\item $D(l\vee l')+D(l\wedge l')=D(l)+D(l')$
\item $D(l)=D(l')\Leftrightarrow l\sim l'$
\item $D(l)\leq D(l')\Leftrightarrow l\leq_{\sim}l'$
\item if $\{l_i\}$ \`{e} is a finite or countable independent system, then
$$
D(\bigvee l_i)=\sum D(l_i)
$$
\item $D(L)$ \`{e} is one of the following sets
\begin{enumerate}
\item $\Delta_n=\{0,\frac{1}{n},\frac{2}{n}, ... \frac{n-1}{n},1\}$
\item $\Delta_{\infty}=[0,1]$
\end{enumerate}
\end{enumerate}
We conclude this section proving some easy facts descending by von
Neumann construction. Let, for the rest of the section, $L$ be a
factorial R-lattice. At first we give the following
\begin{defin} $L$ is said of type $I_n$ if it admits a dimension function $D$
such that $D(L)=\Delta_n$. It is said to be of type $II_1$ if it
admits a dimension function $D$ such that $D(L)=\Delta_{\infty}$.
\end{defin}
\begin{rem}
An easy consequence of the correspondence theorem (see
th.\ref{corrispondenza}) is that the notion of type do not depend by
the choice of the $\sim$dimensionable equivalence relation.
\end{rem}
\begin{rem}
Notation $I_n$ agree with the locution ''of type I''. Indeed, by the
following prop.\ref{dimmin} it follows that a factorial R-lattice is
of type $I_n$ (for some $n$) if and only if it is of type $I$.
\end{rem}
\begin{defin}
A reference on $L$ is given by the choice of a maximal family with
respect to the properties: mutually orthogonal and minimal.
\end{defin}
\begin{prop}\label{dimmin}
Let $L$ of type $I_n$. Then
\begin{enumerate}
\item $l\in Min(L)$ if and only if $D(l)=\frac{1}{n}$
\item if $R$ is a reference on $L$, then $|R|=n$.
\end{enumerate}
\begin{proof}
\begin{enumerate}
\item Let $l\in Min(l)$. If it were $D(l)=0$, it would be $l\sim0$
and thus $l=0$. Instead, if it were $D(l)=m/n$, with $m>1$, one
could consider $l'\in L$ such that $D(l')=1/n$. Whence $0\neq
l'\leq_{\sim}l$ and thus $l$ could not be minimal. Conversely it is
clear using the finiteness property of minimal elements.
\item We suppose for example $|R|=m>n$ (the other case is just the same).
Let $l_1,...l_m\in R$. We obtain the following absurd
$$
1=D(1)\geq D(l_1\vee...\vee
l_m)=\sum_{i=1}^mD(l_i)=\sum_{i=1}^m1/n=m/n>1
$$
\end{enumerate}
\end{proof}
\end{prop}
Now we can prove th.\ref{dimensione}.
\begin{proof}
\begin{enumerate}
\item Using Prop.\ref{dimmin}, we have $|I|,|J|\leq n$.
\item If $|I|=|J|$ then $l=\bigvee_il_i\sim\bigvee_jl'_j=l'$ since
$\sim$ is regular. Conversely, we assume for example that
$|I|\leq|J|$, then there exists $J'\subseteq J$ such that $|I|=|J'|$
and thus, by the regularity of $\sim$, we have $\bigvee_{i\in
I}l_i\sim\bigvee_{j\in J'}l_j'=l_1'$. Whence $l_1'\leq l'$ and
$l_1'\sim l'$. By the finiteness property it follows that $l_1'=l'$
and thus, since $l_j'$ are mutually orthogonal, $|J|=|J'|$.
\end{enumerate}
\end{proof}
\section{Factorial $W^*$-lattices of type $I_n$}
Since we want axiomatize the projection lattice of $M_n(\mathbb C)$,
we are interested in the case in which the lattice contains infinite
non-countable elements. This property does not hold in generale,
since one can easily construct finite factorial lattices for example
of type $I_2$: $L=\{0,x,1-x,y,1-y,1\}$ in which $x,1-x,1-y,y$ have
dimension $1/2$. Thus they are minimals and consequently $x$ (resp.
$y$) commutes only with $0,x,1-x,1$ (resp. $0,y,1-y,1$). Whence $L$
is factorial.
\begin{defin}
A factorial $W^*$-lattice of type $I_n$ is a factorial R-lattice of
type $I_n$ that contains infinite non-countable elements.
\end{defin}
\section{Geometry of minimal elements}\label{minimali}
Let $L$ be a factorial $W^*$-lattice of type $I_n$. We start proving
that also $Min(L)$ has the same cardinality of $L$.
\begin{prop}\label{contmini}
$$
|Min(L)|=|\mathbb R|
$$
\begin{proof}
It follows by Cor.\ref{minimale} and Th.\ref{dimensione} that each
element of $L$ is sup of a finite family of minimal elements. Thus
$L$ has the same cardinality of the set of the finite subsets of
$Min(L)$, which has the same cardinality of $Min(L)$.
\end{proof}
\end{prop}
\begin{defin}
Let $A$ be a set, $a_0\in A$ and $A'\subseteq A\times A$. We say
that $A'$ is $a_0$-right-separated  if $(A\times\{a_0\})\cap
A'=\emptyset$.
\end{defin}
\begin{lem}\label{insiemi}
Let $A, B$ be two equipotent sets, $A'\subseteq A$, $B'\subseteq B$
be also equipotent, $A''\subseteq A\times A$ $a_0$-right-separated
and equipotent to the $b_0$-right-separated set $B''\subseteq
B\times B$. Then there exists a bijection $\Psi:A\rightarrow B$ such
that
\begin{enumerate}
\item $\Psi|_{A'}:A'\rightarrow B'$ is a bijection
\item $(a_1,a_2)\in A''\Leftrightarrow(\Psi(a_1),\Psi(a_2))\in B''$
\end{enumerate}
\begin{proof}
Let us consider the following subsets of $A\times A$:
$A_1=A'\times\{a_0\}, A_2=(A\backslash A')\times\{a_0\}, A_3=A'',
A_4=(A\times A)\backslash(A''\cup A\times\{a_0\})$. They are a
partition of $A\times A$ and they are respectively equipotent to
$B_1=B'\times\{b_0\}, B_2=(B\backslash B')\times\{b_0\}, B_3=B'',
B_4=(B\times B)\backslash(B''\cup B\times\{b_0\})$, which is a
partition of $B\times B$. Thus there exist bijections
$\psi_1:A_1\rightarrow B_1$, $\psi_2:A_2\rightarrow B_2$,
$\psi_3:A_3\rightarrow B_3$ and $\psi_4:A_4\rightarrow B_4$. Now we
can join all these maps to obtain a bijection $\Psi_1:A\times
A\rightarrow B\times B$, which, restricted to $A\times\{a_0\}$,
gives a bijection between it (and thus between $A$) and
$B\times\{b_0\}$ (and thus with $B$) that satisfies the required
properties.
\end{proof}
\end{lem}
It follows an unexpected and decisive result.
\begin{defin}
Let $I$ be a set and $l^2(I)$ the standard Hilbert space on $I$. Let
$l^2(I)_1$ be the unit sphere of $l^2(I)$, that is
$$
l^2(I)_1=\{x\in l^2(I):||x||=1\}
$$
We define the following equivalence relation on $l^2(I)_1$: $x E
y\Leftrightarrow \exists\theta\in[0,2\pi):x=e^{i\theta}y$, in which
transitivity follows by summing $\theta$ modulo $2\pi$.
\end{defin}
\begin{teo}\label{elledue}
Let $L$ be a factorial $W^*$-lattice of type $I_n$ on which a
reference $R=\{r_1,...r_n\}$ is fixed. Then there exists a bijection
$\Psi:Min(L)\rightarrow l^2(R)_1/E$ such that
\begin{enumerate}
\item $\Psi(r_i)$ is the i-th element of the canonical orthonormal basis of the Hilbert space $l^2(R)$.
\item $m\perp n\Leftrightarrow \Psi(m)\perp\Psi(n)$ in the Hilbert space $l^2(R)$.
\end{enumerate}
\begin{proof}
Let $A=Min(L)\cup\{1\}\subseteq L, A'=R, A''=\{(m,n)\in A\times
A:m\perp n\}, B=l^2(R)_1/E\cup\{1\}\subseteq B(l^2(R)),
B'=$''canonical orthonormal basis of $l^2(R)$'', $B''=\{(x,y)\in
B\times B:x\perp y\}$. It is enough to prove that those sets satisfy
the hypothesis of lemma \ref{insiemi} and that $\Psi(1)=1$. In this
case we can restrict $\Psi$ and obtain a bijection
$\Psi:A\rightarrow B$ without modifying its properties. Certainly
$A''$ and $B''$ are $1$-right-separated. Now, since $R$ is finite,
then $l^2(R)_1$ has the same cardinality of $\mathbb R$, that is the
same cardinality of $A$ (by Prop.\ref{contmini}). Moreover $A'$ is
equipotent to $B'$ since each basis of $l^2(R)$ has the same
cardinality of $R$ and thus the same cardinality of $A'$. Lastly
$A''$ is equipotent to $A$ and $B''$ is equipotent to $B$ and thus
$A$ is equipotent to $B$. It remains to prove that $\Psi(1)=1$: it
is sufficient to observe that in the proof of lemma \ref{insiemi}
one can always choose $\psi_1$ such that
$\psi_1(a_0,a_0)=(b_0,b_0)$.
\end{proof}
\end{teo}
After this fundamental result we prove some easy lemmas on the
behavior of the minimal elements. Also these lemmas will be crucial
in the proof of the correspondence theorem.
\begin{lem}\label{zero}
For each $l\in L$, one has $\bigvee_{m\in Min(l)}m=l$.
\begin{proof}
Since $m\leq l$ for each $m\in Min(l)$, one certainly has
$\bigvee_{m\in Min(l)}m\leq l$. Conversely, if $l-\bigvee_{m\in
Min(l)}m\neq0$, then there exists $m'\in Min(l-\bigvee_{m\in
Min(l)}m)$. This is a contradiction, because we have already get all
the minimal elements dominated by $l$.
\end{proof}
\end{lem}
\begin{lem}\label{primo}
Let $l,l'\in L$. Then $l\leq l'$ if and only if $Min(l)\subseteq
Min(l')$.
\begin{proof}
If $l\leq l'$ and $m\in Min(l)$ one certainly has $m\in Min(l')$.
Conversely, using lemma \ref{zero}, one has
$$
l=\bigvee_{m\in Min(l)}m\leq\bigvee_{m\in Min(l')}m=l'
$$
\end{proof}
\end{lem}
\begin{lem}\label{secondo}
Let $l,l'\in L$. Then $\perp(l\vee
l')\subseteq\perp(l)\cap\perp(l')$.
\begin{proof}
Let $x\in\perp(l\vee l')$, then, by definition, $x\leq1-(l\vee l')$.
Now, we know that $l,l'\leq l\vee l'$ and thus $1-l\vee
l'\leq1-l,1-l'$. So $x\leq1-l,1-l'$ and $x\in\perp(l)\cap\perp(l')$.
\end{proof}
\end{lem}
\begin{lem}\label{terzo}
Let $l,l'\in L$. Then
$$
Min(l\vee l')\cap\perp(l)\cap\perp(l')=\emptyset
$$
\begin{proof}
Let $x\in\perp(l)\cap\perp(l')$, then $x\leq1-l,1-l'$. Consequently
$x\leq(1-l)\wedge(1-l')=1-(l\vee l')$. So $x\in\perp(l\vee l')$ and
thus it can not be dominated by $l\vee l'$.
\end{proof}
\end{lem}
\begin{lem}\label{quarto}
Let $m,n\in Min(L)$. Then $m\sim n$.
\begin{proof}
It must be $m\leq_{\sim}n$ or $n\leq_{\sim}m$. Since each non zero
element is not equivalent to 0 and since $m,n$ are minimal, the it
must be exactly $m\sim n$.
\end{proof}
\end{lem}
\begin{lem}\label{quinto}
For each $l'\leq l$, one has $Min(l-l')=Min(l)\cap\perp(l')$.
\begin{proof}
If $x\in Min(l-l')$ one certainly has $x\in Min(l)$. Moreover $x\leq
l-l'\leq1-l'$ and thus $x\in\perp(l')$. Conversely, if $x\in
Min(l)\cap\perp(l')$, then $x\leq1-l'$ is minimal. But $x\leq l$ and
thus $x\leq(1-l')\wedge l=l-l'$ is minimal.
\end{proof}
\end{lem}
\section{The correspondence theorem}
In this section we finally prove that the concept of factorial
$W^*$-lattice of type $I_n$ axiomatizes the projection lattice of
$M_n(\mathbb C)$.
\begin{nt}
$P_n$ will denote the projection lattice of $B(M_n(\mathbb C))$.
\end{nt}
\begin{teo}\label{corrispondenza}
The map $M_n(\mathbb C)\rightarrow P_n$ sends type $I_n$ factors in
factorial $W^*$-lattices of type $I_n$. Conversely, if $L$ is a
factorial $W^*$-lattice of type $I_n$, then $L\cong P_n$.
\begin{proof}
Of course $P_n$ is a factorial $W^*$-lattice of type $I_n$, using
the normalized trace as dimension. Conversely, let $\{l_i\}_{i=1}^n$
be a reference on $L$ and we set $H=l^2(I)$. We map $l_i$ in the
projection $e_{l_i}$ of $H$ onto the $i$-th addend of the direct sum
$H=\mathbb C^n$. Now, if $l\in L$ is minimal, we map it in the
projection of $H$ onto $\Psi(l)\mathbb C$. Now, if $l\in L$, we can
write $l=\bigvee l_j$ (where the $l_j$ are minimal and mutually
orthogonal). Since all the $\Psi(l_j)$ are still minimal and
mutually orthogonal in $M_n(\mathbb C)$, we can map $l$ in the
projection $e_l$ of $H$ onto $\bigoplus_{j\in J}\Psi(l_j)\mathbb C$.
In this way, we have constructed a map from $L$ into $P_n$. Now, if
$e\in P_n$, we write $eH$ as direct sum of one-dimensional subspaces
and then observe that each one-dimensional subspace is of the form
$a\mathbb C$ (for a unique $a\in l^2(I)_1/E$). So, recalling that
$\Psi$ is a bijection, we have that $eH$ comes from
$\bigvee\Psi^{-1}(a)$, where $a$ runs over a decomposition of $eH$.
Thus we have a bijection between $L$ and $P_n$. Now we have to prove
that they are the same lattice structure. Let $\Psi:L\rightarrow
P(B(H))$ be the bijection we have already constructed. We have to
prove the following properties:
\begin{enumerate}
\item $\Psi$ preserves the lattice structure. That is
\begin{enumerate}
\item $\Psi(l)\leq\Psi(l')\Leftrightarrow l\leq l'$
\item $\Psi(l\wedge l')=\Psi(l)\wedge\Psi(l')$
\item $\Psi(l\vee l')=\Psi(l)\vee\Psi(l')$
\end{enumerate}
\item $\Psi$ preserves the orthocomplementation. That is
$$
\Psi(1-l)=\Psi(1)-\Psi(l)
$$
\item $\Psi$ preserves the regular equivalence relation. That is
$$
\Psi(l)\sim\Psi(l')\Leftrightarrow l\sim l'
$$
In this case $\Psi$ would automatically preserve the dimension.
\end{enumerate}
We recall that the map $\Psi|_{Min(L)}$ is a bijection between
$Min(L)$ and $Min(P_n)$ (by construction).
\begin{enumerate}
\item
\begin{enumerate}
\item Let $l\leq l'$ and $\Psi(m)\in Min(\Psi(l))$. Then $m\in
Min(l)\subseteq Min(l')$ and thus $\Psi(m)\in Min(\Psi(l')$.
Consequently, we have $Min(\Psi(l))\subseteq Min(\Psi(l'))$ and, by
lemma \ref{primo}, $\Psi(l)\leq\Psi(l')$. Conversely, if
$\Psi(l)\leq\Psi(l')$, then $Min(\Psi(l))\subseteq Min(\Psi(l'))$.
Thus if $m\in Min(l)$, then $\psi(m)\in Min(\Psi(l))$ and
$\psi(m)\in Min(\Psi(l'))$. Consequently $m\in Min(l')$. So we have
obtained $Min(l)\subseteq Min(l')$ and then $l\leq l'$ (always using
lemma \ref{primo}).
\item Let $\Psi(m)\in Min(P_n)$. We have
$$
\Psi(m)\leq\Psi(l\wedge l')\Leftrightarrow m\leq l\wedge
l'\Leftrightarrow m\leq
l,l'\Leftrightarrow\Psi(m)\leq\Psi(l),\Psi(l')\Leftrightarrow
$$
$$
\Leftrightarrow\Psi(m)\leq\Psi(l)\wedge\Psi(l')
$$
So $\Psi(l\wedge l')$ and $\Psi(l)\wedge\Psi(l')$ have the same
minimal elements. Whence they are equals, using lemma \ref{primo}
again.
\item Certainly $l,l'\leq l\vee l'$. Thus
$\Psi(l),\Psi(l')\leq\Psi(l\vee l')$ and consequently
$\Psi(l)\vee\Psi(l')\leq\Psi(l\vee l')$. Conversely, let $\Psi(m)\in
Min(\Psi(l\vee l')-\Psi(l)\vee\Psi(l'))$, then
$$
\Psi(m)\in Min(\Psi(l\vee
l'))\cap\perp(\Psi(l)\vee\Psi(l'))\subseteq
$$
using lemma \ref{secondo}
$$
\subseteq Min(\Psi(l\vee l'))\cap\perp(\Psi(l))\cap\perp(\Psi(l'))
$$
Whence
$$
m\in Min(l\vee
l')\cap\perp(l)\cap\perp(l')=\emptyset\,\,\,\,(by\,\,\ref{terzo})
$$
That is absurd. Consequenlty $\Psi(l\vee l')-\Psi(l)\vee\Psi(l'))$
does not contain minimal elements and thus it must be zero.
\end{enumerate}
\item One has
$$
\Psi(m)\in Min(\Psi(1-l))\Leftrightarrow m\in
Min(1-l)\Leftrightarrow m\in Min(1)\cap\perp(l)\Leftrightarrow
$$
$$
\Leftrightarrow\Psi(m)\in
Min(\Psi(1))\cap\perp(\Psi(l))\Leftrightarrow\Psi(m)\in
Min(\Psi(1)-\Psi(l))
$$
Thus $\Psi(1-l)$ and $\Psi(1)-\Psi(l)$ have the same minimal
elements and thus they are equals.
\item It remains to prove that $\Psi$ preserves the regular equivalence relation.
Let $\sim_L$ be the equivalence relation on $L$ and let $\sim_{P_n}$
be the equivalence relation on $P_n$. We have to prove that
$l\sim_Ll'\Leftrightarrow\Psi(l)\sim_{P_n}\Psi(l')$. Let
$l\sim_Ll'$, we write $l=\bigvee_{i\in I}l_i$ e $\l'=\bigvee_{i\in
I}l_i'$, where $\{l_i\}$ and $\{l_i'\}$ are two maximal families of
mutually orthogonal and minimal elements dominated respectively by
$l$ and $l'$. Since $l_i$ and $l_i'$ are minimal, it must be
$l_i\sim_Ll_i'$ and thus, since $\Psi$ preserves minimality,
$\Psi(l_i)\sim_{P_n}\Psi(l_i')$ for each $i$. Now, since $\Psi$
preserves minimality, orthogonality and order, $\{\Psi(l_i)\}$ and
$\{\Psi(l_i')\}$ are two maximal families of minimal projections of
$M_n(\mathbb C)$ dominated respectively by $\Psi(l)$ and $\Psi(l')$.
This means that $\Psi(l)$ and $\Psi(l')$ project onto two subspaces
of the same dimension and thus they are equivalent. Conversely, we
assume that $\Psi(l)$ and $\Psi(l')$ are equivalent and thus that
they project onto subspaces of the same dimension. We decompose
these two subspaces in one-dimensional and mutually orthogonal
subspaces and so we write $\Psi(l)=\bigvee_{i\in I}\Psi(l_i)$ and
$\Psi(l')=\bigvee_{i\in I}\Psi(l_i')$, where $\{l_i\}$ and
$\{l_i'\}$ are two maximal families of minimal and mutually
orthogonal elements dominated respectively by $l$ and $l'$. Thus
$l_i\sim l_i'$ for each $i$ and consequently $l=\bigvee
l_i\sim_L\bigvee l_i'=l'$.
\end{enumerate}
\end{proof}
\end{teo}
\begin{rem}
Some easy consequences of the correspondence theorem:
\begin{enumerate}
\item In a factorial $W^*$-lattice of type I the parallelogram rule
is automatically verified.
\item In a factorial $W^*$-lattice of type I there is only one
dimensionable regular equivalence relation, the perspective one.
\end{enumerate}
\end{rem}
\begin{rem}
This remark is suggested by F. Radulescu (\cite{Ra}).\\
Let $L_n$ be a $W^*$-lattice of type $I_n$ and $P_n$ the projections
lattice of $M_n(\mathbb C)$. The correspondence theorem guarantees
$L_n\cong P_n$. A first consequence is the not obvious inclusion
$L_n\hookrightarrow L_{n+1}$. Not obvious because it is not
sufficient to chose those elements in $L_{n+1}$ whose dimension is
less then $\frac{n}{n+1}$: this is not a lattice! Actually this
inclusion does not depend just by the correspondence theorem, but by
the existence of references: choose in $L_{n+1}$ a reference
$R=\{m_1,...m_{n+1}\}$, the sublattice generated by $\{m_1,...m_n\}$
is isomorphic to $L_n$. This argument shows also that it is possible
to choose references on $L_n$ and $P_n$ such that each square of the
following diagram commute.
$$
\xymatrix{L_1\ar[r]^{\xi_1}\ar[d]& P_1\ar[d] \\
L_{2}\ar[r]^{\xi_2}\ar[d]&P_2\ar[d]\\
...\ar[d]&...\ar[d]\\
L_n\ar[r]^{\xi_n}\ar[d]&P_n\ar[d]\\
L_{n+1}\ar[r]^{\xi_{n+1}}\ar[d]&P_{n+1}\ar[d]\\
...&... }
$$
$\xi_n$ being the isomorphisms induced by the choice of increasing
references. Now, since the hyperfinite factor of type $II_1$ is the
inductive limit of $M_n(\mathbb C)$ with normalized trace, we can
find a factorial $W^*$-lattice of type $II_1$ as inductive limit of
factorial $W^*$-lattices of type $I_n$. Moreover, this inductive
limit coincides with the direct limit (with respect the canonical
inclusion) completed by cuts (Dedekind-MacNeille completing, see
\cite{MacN}).
\end{rem}
\begin{prob}
Does the concept of $W^*$-lattice of type I axiomatize the
projections lattice of a type I finite von Neumann algebra one?
\end{prob}
\section{A straightening theorem?}
In this section we want to describe a possible improvement of the
correspondence theorem. Indeed we believe it is possible to remove
hypothesis on the existence of an orthogonal complement. This hope
can be seem exaggerate and thus vain, but actually there are many
cases in which something similar happens. The main one, since it
looks like the our case, is the following: a complex linear space
$V$ of dimension $n$ is isomorphic to $\mathbb C^n$, which has a
notion of orthogonality that we can transport on $V$. On the other
hand, if we analyze our construction we notice that the notion of
orthogonality is been used only in two step: the first one is to
define central elements and then the abelian ones. Nevertheless the
notion of factoriality, following von Neumann, does not really
depend by the existence of central elements, but it mainly depends
by the uniqueness of the inverse (see Def.\ref{fattore}); moreover,
abelian elements allow only to prove the existence of minimal
elements. Thus we can require the existence of a minimal element as
an axiom. The second step in which we used the orthogonality is to
define references. This is just what happens in the case of a
generic linear space of dimension $n$ with respect to $\mathbb C^n$:
in a similar way, we believe it is possible to start from an affine
reference (in the sense of the following definition) and to find a
lattice isomorphism with a factorial $W^*$-lattice of type $I_n$,
$n$ being the cardinality (invariant!) of the affine reference.
\begin{defin}\label{affine}
An affine reference is a family of minimal elements $R=\{l_i\}$ such
that
\begin{enumerate}
\item $\bigvee l_i=1$
\item Each subset of $R$ does not verify the first condition.
\end{enumerate}
\end{defin}
\begin{defin}
An irreducible and continuous lattice with minimal elements is an
irreducible complete lattice such that $|L|=|\mathbb R|$ and
$Min(L)\neq\emptyset$.
\end{defin}
Let $L$ be an irreducible and continuous lattice with minimal
elements.
\begin{rem}
The second condition of Def.\ref{affine} implies that $\{l_i\}$ are
completely independents. Indeed, if not, let $J$ and $i_0$ such that
$l_{i_0}\leq\bigvee_{j\in J}l_j$ and $i_0\notin J$. Then
$\bigvee_{i\in I\setminus\{i_0\}}l_i$ is still equal to 1.
Consequently (in the case of irreducible and complete lattices) we
can still apply th.\ref{dimensione} and obtain that two affine
references have the same finite cardinality.
\end{rem}
\begin{prop}\label{minimal2}
Each element of $L$ is sup of a family of minimal elements.
\begin{proof}
Applying Prop. \ref{mini} (choosing $l$ minimal instead abelian!) we
have that each element $l\in L$ contains minimal elements. If
$\bigvee_{l'\in Min(l)}l'<l$ we apply a classical lemma of J. von
Neumann: if $l'\leq l$, there exists $l''$ (independent to $l'$!)
such that $l'\vee l''=l$. Take a minimal element contained in $l'$
and find the absurd.
\end{proof}
\end{prop}
\begin{con}{\bf (Straightening theorem?)}\label{dritto}
Let $L$ be an irreducible and continuous lattice with minimal
elements and let $n$ be the cardinality of one of its affine
reference. There exists an involutory anti-automorphism
$\perp:L\rightarrow L$ such that $(L,\leq,\perp)$ is a factorial
lattice of type $I_n$.
\end{con}
\begin{rem}
This should be the best theorem we can find, because it is quite
simple to prove that the axioms of irreducible and continuous
lattice with minimal elements are independent: we can find examples
of lattices satisfying all of them except one.
\end{rem}
\section{The separable case}
The main obstacle that we can find in order to extend the
correspondence theorem to the separable case is the absence of the
modular property. But it is been used only to prove that if $l'\leq
l$, one can define $l-l'$. And this last property is been used only
to prove that if there exists a minimal element, then each element
is sup of a family of minimal elements. We think that there are no
other way that assume this property as an axiom.
\begin{defin}\label{omega}
A factorial $W^*$-lattice of type $I_{\omega}$ is a factorial
$W^*$-lattice in which
\begin{enumerate}
\item Each element is sup of a family of minimal elements.
\item 1 is sup of a countable family of minimal elements.
\end{enumerate}
\end{defin}
\begin{teo}
There exists only one (up to lattice isomorphism) factorial
$W^*$-lattice of type $I_{\omega}$ and it is the projection lattice
of $B(H)$, where $H$ is a separable Hilbert space.
\begin{proof}
The proof is the same as in the finite dimensional case. Indeed
results in Sect.\ref{minimali} still hold and we can still make the
construction of the proof of the correspondence theorem thanks to
the two properties 1. and 2. of Def.\ref{omega}.
\end{proof}
\end{teo}
We can also conjectured the straightening theorem in the separable
case.
\begin{con}
Let $L$ be an irreducible, complete and continuous lattice such that
\begin{enumerate}
\item Each element is sup of a family of minimal elements.
\item 1 is sup of a countable family of minimal elements.
\end{enumerate}
Then there exists an involutory anti-automorphism
$\perp:L\rightarrow L$ such that $(L,\leq,\perp)$ is a factorial
lattice of type $I_{\omega}$.
\end{con}
We want to thank Prof. F. Radulescu and Prof. L. Zsido for the
useful discussion.

VALERIO CAPRARO, Universit\`{a} degli Studi di Roma "Tor Vergata".
Email: capraro@mat.uniroma2.it

\end{document}